\pdfoutput=1
\documentclass{mathprint} %
\usepackage{rutarmacros}
\usepackage{macros}

\title[Regularity of non-autonomous self-similar sets]{Regularity of non-autonomous self-similar sets}
\author
  {Antti Käenmäki}
  {University of Eastern Finland, 
   Department of Physics and Mathematics,
   P.O.\ Box 111,
   FI-80101 Joensuu,
   Finland}
  {antti@kaenmaki.net}

\author
  {Alex Rutar}
  {Department of Mathematics and Statistics,
      University of Jyväskylä,
      P.O.\ Box 35 (MaD),
      FI-40014 University of Jyväskylä,
  Finland}
  {alex@rutar.org}

\begin{document}
\begin{abstract}
    Non-autonomous self-similar sets are a family of compact sets which are, in some sense, highly homogeneous in space but highly inhomogeneous in scale.
    The main purpose of this note is to clarify various regularity properties and separation conditions relevant for the fine local scaling properties of these sets.
    A simple application of our results is a precise formula for the Assouad dimension of non-autonomous self-similar sets in $\mathbb{R}^d$ satisfying a certain ``bounded neighbourhood'' condition, which generalizes earlier work of Li--Li--Miao--Xi and Olson--Robinson--Sharples.
    We also see that the bounded neighbourhood assumption is, in few different senses, as general as possible.
\end{abstract}

\section{Introduction}
In the study of non-smooth and highly irregular sets, one encounters sets with a wide variety of structure.

On one hand, one has sets with a substantial amount of structure: perhaps the most canonical example of such a set is the \emph{middle-thirds Cantor set}.
Many generalizations exist, such as conformal or expanding repellers, sets occurring in complex dynamics, and random sets and often those occurring in models of physical phenomena.

At the other extreme, we have the class of general compact, or even analytic, sets.
For such sets, one is restricted to only making the most general observations because of an abundance of counterexamples.
Moreover, such general observations are often very challenging to make and depend on deep properties of the ambient space.

It turns out that a useful intermediate family of sets are those with inhomogeneity in scale, but homogeneity in space.
Such sets are commonly referred to in the literature as \emph{Moran sets} or attractors of \emph{non-autonomous iterated function systems}.
Beyond being of interest in their own right, this family of sets has also played an important role of being sufficiently unstructured to be applicable in a wide variety of scenarios, while still being sufficiently structured to make the analysis of their properties convenient.
Some examples of applications in which this class of sets is highly useful include the construction of large non-ergodic measures \cite{zbl:1184.28009}, the dimension theory of non-uniform sets arising in complex dynamics \cite{zbl:1500.37035}, and classification results for dimensions \cite{arxiv:2206.06921,zbl:1509.28005}.
For the authors, this paper was motivated by applications to slices of self-affine sets \cite{arxiv:2410.19404}.

In this paper, we study non-autonomous self-similar sets in their own right.
Our focus is on coarse notions of dimension: in particular, we focus on the \emph{Assouad dimension}.
The Assouad dimension was first introduced in the study of bi-Lipschitz embeddings of metric spaces \cite{zbl:0396.46035}.
It also previously appeared implicitly in Furstenberg's work on dynamics on fractals as the \emph{star dimension} \cite{zbl:0208.32203,zbl:1154.37322}.
Especially in recent years, a large amount of attention has been paid to the Assouad dimension.
We refer the reader to the books \cite{zbl:1201.30002,zbl:1467.28001,zbl:1222.37004} and the many references therein for an introduction to this specific area from a variety of perspectives.

Our primary goal is to organize the study of non-autonomous self-similar sets around a few natural definitions, which are relevant in a much broader context.
We hope that with the definitions in mind, for the experienced reader, the proofs will follow naturally.
A consequence of this (re)organization is that we obtain essentially optimal generalizations of existing results in the literature concerning the Assouad dimension of non-autonomous self-similar sets, and moreover almost for free.

Let us begin by introducing our setting more precisely.

\subsection{Non-autonomous self-similar sets}
The framework that we find most natural is that of the \emph{non-autonomous conformal iterated function system}, which was first introduced in \cite{zbl:1404.28017}.
In that paper, under certain regularity assumptions, Rempe-Gillen \& Urbański prove that the Hausdorff and box dimensions are equal and given by the zero of a certain pressure function.

We specialize slightly and consider only self-similar systems.
Generalization to conformal systems follows without difficulty.
For each $n\in\N$, let $\mathcal{J}_n$ be a finite index set with $\#\mathcal{J}_n\geq 2$, and let $\Phi_n=\{S_{n,j}\}_{j\in\mathcal{J}_n}$ be a family of similarity maps $S_{n,j}\colon\R^d\to\R^d$ of the form
\begin{equation*}
    S_{n,j}(\bm{x}) = r_{n,j} O_{n,j}\bm{x} + \bm{d}_{n,j}
\end{equation*}
where $r_{n,j}\in(0,1)$ and $O_{n,j}$ is an orthogonal matrix.
To avoid degenerate situations, we assume that associated with the sequence $(\Phi_n)_{n=1}^\infty$ is an invariant compact set $X\subset\R^d$ (that is $S_{n,j}(X)\subset X$ for all $n\in\N$ and $j\in\mathcal{J}_n$) and moreover that
\begin{equation}\label{e:rn-lim}
    \lim_{n\to\infty}\sup\{r_{1,j_1}\cdots r_{n,j_n}:j_i\in\mathcal{J}_i\text{ for each }i=1,\ldots,n\}=0.
\end{equation}
Under these assumptions, associated with the sequence $(\Phi_n)_{n=1}^\infty$ is a \defn{limit set}
\begin{equation*}
    K=\bigcap_{n=1}^\infty\bigcup_{(j_1,\ldots,j_n)\in\mathcal{J}_1\times\cdots\times\mathcal{J}_n} S_{1,j_1}\circ\cdots\circ S_{n,j_n}(X).
\end{equation*}
Since $X$ is compact and invariant under any map $S_{n,j}$ with $j\in\mathcal{J}_n$, finiteness of each $\mathcal{J}_n$ implies that $K$ is the intersection of a nested sequence of compact sets and therefore non-empty and compact.
Under these assumptions, the sequence $(\Phi_n)_{n=1}^\infty$ is called a \emph{non-autonomous iterated function system (IFS)} and the limit set $K$ is called the \emph{non-autonomous self-similar set}.
We say that $(\Phi_n)_{n=1}^\infty$ is \emph{autonomous} if the $\Phi_n$ do not depend on the index set $n$, and \emph{homogeneous} if for each $n\in\N$, there is a number $r_n$ so that $r_{n,j}=r_n$ for all $j\in\mathcal{J}_n$.

We emphasize, at this point, that we require no other technical assumptions concerning the set $X$ or the contraction ratios $r_{n,j}$ as can be found in \cite{zbl:1404.28017}.
We refer the reader to \cite[§2]{zbl:1404.28017} for more detail on this construction in the conformal setting.

Closely related to our setup is the slightly more general \emph{Moran set construction}; see, for instance, \cite{zbl:1148.28007,zbl:0881.28003}.
In the Moran set construction, the placement of the geometric cylinders is allowed to vary across a given level while still respecting separation conditions.
Our results also hold for this more general class of sets, but for clarity of notation, we prefer the framework of non-autonomous iterated function systems.

\subsection{Symbolic formulas for dimension}\label{ss:symb}
Let $n\in\N$ and $m\in\N$ be arbitrary.
We define the \emph{pressure} $\phi_{n,m}\colon \R\to\R$ by
\begin{equation*}
    \phi_{n,m}(t)= \frac{1}{m} \log\sum_{j_1\in\mathcal{J}_{n}}\cdots\sum_{j_m\in\mathcal{J}_{n+m-1}}\prod_{k=0}^{m-1} r_{n+k,j_{n+k}}^t.
\end{equation*}
It is easy to check that $\phi_{n,m}$ is convex, differentiable, and strictly decreasing with unique zero $\theta(n,m)\geq 0$.
Of course, $\theta(n,m)$ is precisely the similarity dimension of the IFS
\begin{equation*}
    \Phi_{n+1}\circ\cdots\circ\Phi_{n+m}=\{f_1\circ\cdots\circ f_m:f_i\in\Phi_{n+i}\}.
\end{equation*}
A natural question is the following: \emph{when are the dimensions of the limit set $K$ fully determined by the numbers $\theta(n,m)$?}
In an ideal world, some notions of dimensions might be given by the formulas
\begin{itemize}[nl]
    \item $\dimuB K = \limsup_{m\to\infty}\theta(1,m)$.
    \item $\dimH K = \liminf_{m\to\infty}\theta(1,m)$.
    \item $\dimA K = \limsup_{m\to\infty}\sup_{n\in\N}\theta(n,m)$.
\end{itemize}
Of course, without more assumptions on the IFS $(\Phi_n)_{n=1}^\infty$, there is no chance that any of these formulas can be true.
Even in the presence of separation, by taking the $\mathcal{J}_n$ to be very large, we may essentially approximate any compact subset of $\R^d$ at arbitrary precision and at infinitely many scales.
With this in mind, assumptions on non-autonomous self-similar sets might be grouped into the following categories:
\begin{itemize}[nl]
    \item\label{im:sep}\emph{Separation}: to ensure that the geometric cylinders $S_{1,i_1}\circ\cdots\circ S_{n,i_n}(X)$ do not overlap too much.
    \item\label{im:growth}\emph{Growth rate control}: to ensure that the sizes of the index sets $\mathcal{J}_n$ do not grow too quickly, or at all.
    \item\label{im:scaling}\emph{Bounds on scaling rate}: to ensure that the contraction ratios $r_{n, i_n}$ are not too close to 1, or not too close to 0.
\end{itemize}
For instance, in \cite{zbl:1404.28017}, such a formula is established for Hausdorff dimension under the assumptions (1) that the open set condition holds with respect to a relatively ``nice'' open set; (2) that the index sets have sub-exponential growth rate; and (3) that the scaling ratios are uniformly bounded away from 1.
Through some examples, they also show that these assumptions are essentially optimal.
For ordinary (autonomous) self-similar sets, of course one is mainly concerned with separation, since in this setting one has essentially optimal growth rate control and bounds on the scaling rate.
Such a formula has also previously been established for the Assouad dimension \cite{zbl:1364.28011,zbl:1371.28024} under the assumptions (1) that the open set condition holds with respect to a relatively ``nice'' open set; (2) that the index sets are uniformly bounded in size; and (3) that the scaling ratios are uniformly bounded away from 0 and 1.

One of the key observations in this paper is that all of the three assumptions (separation, growth rate control, and bounds on scaling rate) can instead be replaced with a single, more general, type of assumption: \emph{controlling overlap counts}.

\subsection{Assouad dimension and the bounded neighbourhood condition}\label{ss:bnc-intro}
In this paper, our focus is on notions of dimension which are more sensitive to local scaling: in particular, we focus on the \emph{Assouad dimension}
\begin{align*}
    \dimA K=\inf\Bigl\{s:\exists C>0\,&\forall 0<r\leq R<1\,\forall x\in K\\*
                                         &N_r(B(x,R)\cap K)\leq C\Bigl(\frac{R}{r}\Bigr)^s\Bigr\},
\end{align*}
where $N_r(A)$ is the least number of closed balls of radius $r>0$ needed to cover a bounded set $A \subset \R^d$.
Since the Assouad dimension is very sensitive to the local scaling of the set, we will, in principle, require stronger assumptions than those that are required for Hausdorff dimension.

First, we need some more notation.
Let $\mathcal{T}=\bigcup_{n=0}^\infty\mathcal{T}_n$ denote the set of all cylinder sets in the infinite product space $\Delta=\prod_{n=1}^\infty\mathcal{J}_n$, where
\begin{equation*}
    \mathcal{T}_n = \left\{[j_1,\ldots,j_n] = \{j_1\}\times\cdots\times\{j_n\}\times\prod_{k=n+1}^\infty \mathcal{J}_k:j_i \in\mathcal{J}_i\quad\text{for}\quad i=1,\ldots,n\right\}.
\end{equation*}
Note that the unique cylinder in $\mathcal{T}_0$ is the set $\Delta$.
Given a cylinder $Q=[j_1,\ldots,j_n]\in\mathcal{T}$, we write $\rho(Q) = r_{1,j_1}\cdots r_{n,j_n}$, and if $Q\neq\Delta$, we also let $\widehat{Q}=[j_1,\ldots,j_{n-1}]$ denote the \emph{parent} of $Q$.
The valuation $\rho$ induces a metric $d$ on $\Delta$, given by $d(x,y)=\inf\{\rho(Q):\{x,y\}\subset Q\}$.

Since the contraction ratios need not be the same, it is also natural to stratify the set $\mathcal{T}$ by size.
For $r>0$, we write
\begin{equation*}
    \mathcal{T}(r)= \bigl\{Q\in\mathcal{T}\setminus\{\Delta\}:\rho(Q)\leq r < \rho(\widehat{Q})\bigr\}.
\end{equation*}
We also denote the natural coding map by $\pi\colon\Delta\to\R^d$, which is defined by
\begin{equation*}
    \{\pi\bigl((i_n)_{n=1}^\infty\bigr)\} = \bigcap_{n=1}^\infty S_{1,i_1}\circ\cdots\circ S_{n,i_n}(X).
\end{equation*}
This function is well-defined by \cref{e:rn-lim}, and it is easy to see that $\pi$ is Lipschitz so that $\pi(\Delta)=K$.
\begin{definition}
    Suppose $(\Phi_n)_{n=1}^\infty$ is a non-autonomous IFS.
    We define the \emph{neighbourhood} at $x\in K$ and $r>0$ by
    \begin{equation*}
        \mathcal{N}(x,r)\coloneqq \{Q\in \mathcal{T}(r):\pi(Q)\cap B(x,r)\neq\varnothing\}.
    \end{equation*}
    We then say that the IFS satisfies the \emph{bounded neighbourhood condition} if
    \begin{equation*}
        \limsup_{r\to 0}\sup_{x\in K}\#\mathcal{N}(x,r)<\infty.
    \end{equation*}
\end{definition}
Actually, the bounded neighbourhood condition is almost equivalent to saying that the coding map $\pi$ is bi-Lipschitz.
The correct generalization of a bi-Lipschitz function to make such a statement precise can be found in \cref{ss:coarse-bi} and, in particular, \cref{c:symb-Assouad}.

Our main result is the following, which generalizes previously known results \cite{zbl:1364.28011, zbl:1371.28024}.
\begin{itheorem}\label{it:dima}
    Suppose $(\Phi_n)_{n=1}^\infty$ is a non-autonomous IFS satisfying the bounded neighbourhood condition and $K$ is the non-autonomous self-similar set.
    Then
    \begin{equation}\label{e:dima-form}
        \dimA K = \lim_{m\to\infty}\sup_{n\in\N}\theta(n,m)=\lim_{m\to\infty}\limsup_{n\to\infty}\theta(n,m)=\inf_{m\in\N}\limsup_{n\to\infty}\theta(n,m)
    \end{equation}
\end{itheorem}
Note that one cannot replace the bounded neighbourhood assumption with any unbounded growth of the local covering numbers; see \cref{ex:unbounded}.

Actually, \cref{it:dima} is a consequence of three simple, but in our opinion conceptually useful, observations.
\begin{enumerate}[nl]
    \item The bounded neighbourhood condition holds if and only if the coding map $\pi$ is \emph{bi-Lipschitz decomposable} (see \cref{d:coarse-bi-lip}), which is a generalization of the bi-Lipschitz assumption which in some sense avoids topological obstructions.
        This enables the reduction to the symbolic space.
    \item The limits in \cref{e:dima-form} exist since the function $\theta$ has certain weak quasiconvexity and semi-continuity properties, which generalize the more usual notion of subadditivity.
        We call this property \emph{submaximality}; see \cref{d:submax} for the definition.
    \item It is sufficient to consider only cylinder sets at a fixed level (as opposed to diameter) as a consequence of a \emph{disc-packing} formulation of Assouad dimension, similar to that of the analogous version for box dimension---see, for example, \cite[§2.6]{zbl:1390.28012}.
\end{enumerate}
Unlike previous work, the main technical difficulty in the proof of \cref{it:dima} is that we only have very weak control on the contraction rates.

We hope that our proof helps to clarify the various assumptions which are relevant to the study of non-autonomous self-similar sets, and to unify the proofs of somewhat weaker versions of this result present in the literature.
\begin{remark}
    In the autonomous case, it is a consequence of \cite[Theorem~1.3]{zbl:1317.28014} that if $\dimH K<1$, then \cref{e:dima-form} holds if and only if the bounded neighbourhood condition holds.
    The condition $\dimH K<1$ is needed to avoid saturation of non-trivial affine subspaces.
\end{remark}

\subsection{Other assumptions on non-autonomous IFSs}
In this section, we discuss the difference between the bounded neighbourhood condition and a much more common assumption: the open set condition.
\begin{definition}\label{d:osc}
    We say that a non-autonomous IFS satisfies the \defn{open set condition} if the invariant compact set $X$ can be chosen to have non-empty interior $U=X^\circ$ so that for each $n\in\N$ and $j\neq j'\in\mathcal{J}_n$, we have $S_{n,j}(U)\cap S_{n,j'}(U)=\varnothing$.
    We moreover say that it satisfies the \emph{cone condition} if
    \begin{equation*}
        \inf_{x\in X}\inf_{r\in(0,1)}r^{-d}\mathcal{L}^d\bigl(B(x,r)\cap U\bigr)>0.
    \end{equation*}
\end{definition}
Here, $\mathcal{L}^d$ denotes the usual $d$-dimensional Lebesgue measure.
The cone condition is typically an implicit global assumption \cite{zbl:1404.28017,zbl:0852.28005} when defining a non-autonomous IFS, but we intentionally mention it separately here to make it clear exactly when it is required.

In order to fully explain the distinction between the bounded neighbourhood condition and the open set condition, we make the following definition.
As explained in \cref{ss:symb}, separation in itself is insufficient to guarantee that the desired symbolic formulas hold.
\begin{definition}\label{d:bbranch}
    Given a non-autonomous IFS $(\Phi_n)_{n=1}^\infty$, we define the \emph{geometric offspring count} for $r>0$ by
    \begin{equation*}
        \beta(r)\coloneqq \sup_{Q\in\mathcal{T}(r)}\#\{Q'\in\mathcal{T}(r):Q'\subset \widehat{Q}\}.
    \end{equation*}
    We then say that a non-autonomous IFS has \defn{bounded branching} if
    \begin{equation*}
        \limsup_{r\to 0}\beta(r)<\infty.
    \end{equation*}
\end{definition}
A reasonable way to think about the set $\{Q'\in\mathcal{T}(r):Q'\subset \widehat{Q}\}$ is the set of ``geometric'' siblings of the cylinder $Q\in\mathcal{T}(r)$, as opposed to the set of children of the parent of $Q$.

A simple consequence of bounded branching is that $\sup_{n\in\N}\#\mathcal{J}_n<\infty$ (see \cref{l:branch-bound}).
If the IFS is homogeneous, then it is also easy to check that these two conditions are equivalent.

It is easy to see that there is a non-autonomous IFS which satisfies the bounded neighbourhood condition but not the open set condition; indeed, this will necessarily happen if there are distinct $(i_1,\ldots,i_n)$ and $(j_1,\ldots,j_k)$ such that
\begin{equation*}
    S_{1,i_1}\circ\cdots\circ S_{n,i_n} =  S_{1,j_1}\circ\cdots\circ S_{k, j_k}.
\end{equation*}
To give a concrete example, one might take $\Phi_1 = \Phi_2 = \{x/3, x/3+2/9, x/3+2/3\}$ and $\Phi_n = \{x/3, x/3 + 2/3\}$ for $n\geq 3$.
If the IFS is autonomous, then it is well-known, see e.g.\ \cite[Theorems 3.5 and 5.5]{zbl:1148.28007}, that the bounded neighbourhood condition and the open set condition are equivalent.
\begin{itheorem}\label{it:reg}
    Let $(\Phi_n)_{n=1}^\infty$ be a non-autonomous IFS satisfying the open set condition.
    Then the bounded neighbourhood condition holds if any of the following conditions are satisfied:
    \begin{enumerate}[nl,a]
        \item The IFS satisfies the cone condition and has bounded branching.
        \item The IFS is homogeneous and has bounded branching.
        \item The IFS has contraction ratios uniformly bounded away from 0: that is, there is a constant $r_{\min}>0$ so that $r_{n,j}\geq r_{\min}$ for all $n\in\N$ and $j\in\mathcal{J}_n$.
    \end{enumerate}
\end{itheorem}
The most important aspect of \cref{it:reg} is the statement.
Having stated it, the proof is relatively standard; certainly, the various techniques used in the proof have appeared in previous papers studying non-autonomous iterated function systems.
The point is that these geometric arguments are only required to establish the bounded neighbourhood condition, which, by \cref{c:symb-Assouad}, is essentially equivalent to the reduction to symbolic space and enables one to ignore the overlapping structure of the IFS.

Finally, to demonstrate necessity of our assumptions, we construct a variety of explicit examples with various irregular properties.
In the first example, we show that bounding $\#\mathcal{J}_n$ in the definition of the bounded branching condition is insufficient: more precisely, for any $\varepsilon>0$, there exists a non-autonomous IFS with limit set $K\subset[0,1]$ satisfying the open set condition with respect to the open set $(0,1)$ such that
\begin{enumerate}[nl,r]
    \item has $\#\mathcal{J}_n=2$ for all $n\in\N$,
    \item has $\dimA K=1$, and
    \item has $0<\theta(n,m)\leq\varepsilon$ for all $n,m\in\N$.
\end{enumerate}
In particular, this IFS must be inhomogeneous, or it would satisfy the bounded neighbourhood condition by \cref{it:reg}.
The construction can be found in \cref{ex:unbounded}.

Next, even under the additional assumption of homogeneity, we show that the Assouad dimension does not only depend on the values $\theta(n,m)$.
More precisely, for any sequence $(k_n)_{n=1}^\infty$ with $\limsup_{n\to\infty}k_n=\infty$ and $0<s\leq t\leq 1$, there exists a non-autonomous IFS with limit set $K\subset[0,1]$ satisfying the open set condition with respect to the open set $(0,1)$ such that
\begin{enumerate}[nl,r]
    \item has $\#\mathcal{J}_n\leq k_n$ for all $n\in\N$,
    \item has $\dimA K=t$, and
    \item has $\theta(n,m)=s$ for all $n,m\in\N$.
\end{enumerate}
This example can be found in \cref{ex:arbitrary-values}.

It is reasonable that both of these examples could be modified for even greater flexibility in specifying the values of $\dimA K$ and $\theta(n,m)$.
On the other hand, it is likely that there are certain restrictions concerning the configurations of dimension.
For instance, it is plausible that if the open set condition and the cone condition hold, then $\dimA K \geq \limsup_{m\to\infty}\sup_{n\in\N}\theta(n,m)$.
We believe that this question, and similar questions concerning more general symbolic inequalities, would be relevant directions for subsequent work.

\section{Regularity of non-autonomous self-similar sets} \label{sec:non-autonomous}\label{sec:reduction-bnc}
In this section, we establish a variety of regularity conditions concerning non-autonomous self-similar iterated function systems.
In particular, we show that the bounded neighbourhood condition allows us to ignore the overlapping structure on non-autonomous self-similar sets.
We also study how the bounded neighbourhood condition is related to the open set condition and prove \cref{it:reg}.

\subsection{Lipschitz decomposable relations}\label{ss:coarse-bi}
We recall from \cref{ss:bnc-intro} that the limit set $K$ comes equipped with a Lipschitz coding map $\pi\colon\Delta\to K$.
It would be convenient if the map $\pi$ would preserve Assouad dimension: however, Lipschitz maps may both decrease and increase Assouad dimension.
On the other hand, bi-Lipschitz maps do preserve Assouad dimension.
However, the coding space $\Delta$ is always totally disconnected, whereas $K$ could be a connected set, so in such a situation the coding map cannot be bi-Lipschitz.

The following notion is a generalization of a bi-Lipschitz function between metric spaces.
First, given a relation $\mathcal{R}\subset X\times Y$, its \emph{domain} is the set
\begin{equation*}
    \dom\mathcal{R}=\bigl\{x\in X:(\{x\}\times Y)\cap\mathcal{R}\neq\varnothing\bigr\}.
\end{equation*}
The \emph{image} of a given set $A \subset X$ under the relation $\mathcal{R}$ is the set
\begin{equation*}
  \mathcal{R}(A)=\bigl\{y\in Y:A\times \{y\}\cap\mathcal{R}\neq\varnothing \bigr\}.
\end{equation*}
We also define the \emph{inverse relation} $\mathcal{R}^{-1}=\{(y,x):(x,y)\in\mathcal{R}\}\subset Y\times X$.
\begin{definition}\label{d:coarse-bi-lip}
    Let $(X,d_1)$ and $(Y,d_2)$ be non-empty metric spaces and let $\mathcal{R}\subset X\times Y$ be a relation.
    We say that $\mathcal{R}$ is \emph{Lipschitz decomposable} if $\dom\mathcal{R}=X$ and there are constants $M\in\N$ and $c> 0$ so that for all $x\in X$ and $r>0$, there are $y_1,\ldots,y_M\in Y$ so that
    \begin{equation*}
        \mathcal{R}(B(x,r))\subset\bigcup_{i=1}^M B(y_i, cr).
    \end{equation*}
    We say that $\mathcal{R}$ is \emph{bi-Lipschitz decomposable} if $\mathcal{R}$ and $\mathcal{R}^{-1}$ are both Lipschitz decomposable.
\end{definition}
If we wish to indicate the dependence on the numbers $M$ and $c$, we will say that a relation is $(M,c)$-Lipschitz decomposable.
The following fact is easy to verify.
\begin{lemma}
    Let $(X,d_1)$ and $(Y,d_2)$ be non-empty metric spaces and let $\mathcal{R}\subset X\times Y$ be a relation.
    Then $\mathcal{R}$ is a $c$-Lipschitz function if and only if $\mathcal{R}$ is $(1,c)$-Lipschitz decomposable.
\end{lemma}

In general, pairs of metric spaces which support bi-Lipschitz decomposable relations will have the same dimensions, for most notions of dimensions defined using covers.
Bi-Lipschitz decomposable maps also preserve (up to a multiplication by uniformly bounded constants) quantities defined in terms of covers, such as Hausdorff measure.
For completeness, we give the elementary proof for Assouad dimension and leave consideration of other dimensions to the motivated reader.
\begin{lemma}\label{l:Assouad-lip-cover}
    Let $(X,d_1)$ and $(Y,d_2)$ be metric spaces.
    Suppose there is a bi-Lipschitz decomposable relation $\mathcal{R}\subset X\times Y$.
    Then $\dimA X=\dimA Y$.
\end{lemma}
\begin{proof}
    Choose $M\in\N$ and $c>0$ such that there exists a $(M,c)$-bi-Lipschitz decomposable relation $\mathcal{R}\subset X\times Y$.
    Without loss of generality, we may assume that $c\geq 1$.

    Let $\varepsilon>0$, $x\in X$, and $0<r\leq R<c^{-1}$ be arbitrary.
    By the assumption, get $y_1,\ldots,y_M\in Y$ so that
    \begin{equation}\label{e:img}
        \mathcal{R}(B(x,R))=\{y\in Y:(B(x,R)\times \{y\})\cap\mathcal{R}\neq\varnothing\}\subset\bigcup_{i=1}^M B(y_i, cR).
    \end{equation}
    For each $i=1,\ldots,M$, by definition of the Assouad dimension, get a family of balls $\{B(y_{i,j},c^{-1}r)\}_{i=1}^{N}$ covering $B(y_i,cR)$ where
    \begin{equation*}
        N\lesssim_\varepsilon\left(\frac{cR}{c^{-1}r}\right)^{\dimA Y+\varepsilon}\lesssim \left(\frac{R}{r}\right)^{\dimA Y+\varepsilon}.
    \end{equation*}
    Then for each $B(y_{i,j},c^{-1} r)$, since $\mathcal{R}^{-1}$ is also Lipschitz decomposable, there are balls $\{B(x_{i,j,k},r)\}_{k=1}^M$ so that
    \begin{equation}\label{e:preimg}
        \{x\in X:(\{x\}\times B(y_{i,j},c^{-1}r))\cap\mathcal{R}\neq\varnothing\}\subset\bigcup_{k=1}^M B(x_{i,j,k}, r).
    \end{equation}
    But \cref{e:img} and \cref{e:preimg} together imply that the family of balls
    \begin{equation*}
        \{B(x_{i,j,k},r):i=1,\ldots,M;\, j=1,\ldots, N;\, k=1,\ldots,M\}
    \end{equation*}
    is a cover for $B(x,R)$ with cardinality $M^2 N\lesssim (R/r)^{\dimA Y+\varepsilon}$.
    Since $\varepsilon>0$ was arbitrary, it follows that $\dimA X\leq\dimA Y$.

    Of course, the identical argument applied to the inverse $\mathcal{R}^{-1}$ implies that $\dimA Y\leq \dimA X$, as required.
\end{proof}
\begin{remark}
    We recall that even Lipschitz functions can both decrease and increase Assouad dimension, so we require the full bi-Lipschitz decomposability assumption to establish either inequality.
\end{remark}
It is now straightforward to obtain our result on the Assouad dimension.
\begin{proposition}\label{c:symb-Assouad}
    Suppose $(\Phi_n)_{n=1}^\infty$ is a non-autonomous IFS with limit set $K$ and associated infinite product space $\Delta$.
    Then $(\Phi_n)_{n=1}^\infty$ satisfies the bounded neighbourhood condition if and only if the coding map $\pi$ is bi-Lipschitz decomposable.
    In particular, if either of these equivalent conditions hold, then $\dimA K=\dimA \Delta$.
\end{proposition}
\begin{proof}
    First, suppose the IFS satisfies the bounded neighbourhood condition.
    Since $\pi\colon\Delta\to K$ is Lipschitz, it is Lipschitz decomposable.
    Next, let $M$ be as in the definition of the bounded neighbourhood condition.
    Suppose $y\in K$ and $0<r<1$ are arbitrary: then
    \begin{equation*}
        \#\{Q\in\mathcal{T}(r):\pi(Q)\cap B(y,r)\neq\varnothing\}\leq M.
    \end{equation*}
    But $\diam Q\leq r$, so in fact $\mathcal{R}^{-1}$ is $(M,1)$-Lipschitz decomposable.

    Conversely, suppose $\pi$ is bi-Lipschitz decomposable and let $r>0$ and $y\in K$ be arbitrary.
    Let $M$ and $c$ be such that $\pi^{-1}$ is Lipschitz decomposable.
    Since $\mathbb{R}^d$ is doubling, there exists a constant $N=N(c)$ such that we may cover any ball $B(y,r)$ by $N$ balls of radius $c^{-1} r$.
    Now fix some ball $B(y,r)$ and cover it with balls $B(y_i, c^{-1} r)$ for $i=1,\ldots, N$.
    For each $i$, since $\pi^{-1}$ is Lipschitz decomposable, there exist points $\{x_{i,1},\ldots,x_{i,M}\}\subset\Delta$ such that
    \begin{equation*}
        \pi^{-1}(B(y_i,c^{-1}r))\subset\bigcup_{j=1}^M B(x_{i,j}, r).
    \end{equation*}
    For each $x_{i,j}$, let $Q_{i,j}\in\mathcal{T}(r)$ be the unique cylinder containing $x_{i,j}$.
    Moreover, since $\diam \widehat{Q_{i,j}} > r$, if $x\notin Q_{i,j}$, then $d(x,x_{i,j})>r$.
    Therefore, $B(x_{i,j},r)=Q_{i,j}$ and
    \begin{equation*}
        K \cap B(y,r)\subset K \cap \bigcup_{i=1}^N  B(y_i,c^{-1}r))\subset \bigcup_{i=1}^N\bigcup_{j=1}^M\pi(Q_{i,j}).
    \end{equation*}
    Since $r>0$ and $y\in K$ were arbitrary, it follows that the IFS satisfies the bounded neighbourhood condition with constant $N\cdot M$.

    In particular, if either of these conditions hold, since bi-Lipschitz decomposable maps preserve Assouad dimension by \cref{l:Assouad-lip-cover}, it follows that $\dimA K=\dimA\Delta$.
\end{proof}

\subsection{Consequences of the bounded neighbourhood condition}
In this section, we discuss some of the implications of the bounded neighbourhood condition.
First, we introduce a local non-degeneracy condition similar to \cref{e:rn-lim}.
\begin{definition}
    We say that the IFS is \emph{locally contracting} if
    \begin{equation*}
        \lim_{m\to\infty}\sup_{n\in\N}\max\{r_{n,j_0}\cdots r_{n+m-1,j_{m-1}}:\,j_i\in\mathcal{J}_{n+i}\text{ for each }i=0,\ldots,m-1\}= 0.
    \end{equation*}
\end{definition}
We have the following consequences of the bounded neighbourhood condition.
\begin{lemma}\label{l:bnc-b}
    Suppose $(\Phi_n)_{n=1}^\infty$ is a non-autonomous IFS satisfying the bounded neighbourhood condition.
    Then it is locally contracting and $\sup_{n\in\N}\#\mathcal{J}_n < \infty$.
\end{lemma}
\begin{proof}
    Throughout, let $M\in\N$ be the constant with respect to which the IFS satisfies the bounded neighbourhood condition.

    First, suppose for contradiction that the IFS is not locally contracting.
    Then there is a $\delta>0$ so that for all $m\in\N$, there is an $n\in\N$ and $(i_0,\ldots,i_{m-1})\in\mathcal{J}_n\times\cdots\times\mathcal{J}_{n+m-1}$ such that
    \begin{equation*}
        r_{n,i_0}\cdots r_{n+m-1,i_{m-1}} \geq \delta.
    \end{equation*}
    Fix $m\in\N$ and the corresponding $n\in\N$ and $(i_0,\ldots,i_{m-1})$ as above; we may moreover assume that $r_{n,i_0}\cdots r_{n+m-1,i_{m-1}}$ is maximal.
    Also, choose some arbitrary initial segment $(j_1,\ldots,j_{n-1})$, let $Q_0 = [j_1,\ldots,j_{n-1}]$, and let $Q = [j_1,\ldots,j_{n-1},i_0,\ldots,i_{m-1}]$.

    Write $r=\rho(Q)$.
    Since $\#\mathcal{J}_k\geq 2$ for all $k\in\N$, for each $\ell=0,\ldots,m-1$, there is a cylinder
    \begin{equation*}
        Q'_\ell = [j_1,\ldots,j_{n-1},i_0,\ldots,i_\ell']
    \end{equation*}
    where $i_\ell' \neq i_\ell$.
    Moreover, by maximality of the choice of $[i_0,\ldots,i_{m-1}]$, each $Q'_\ell$ is either itself an element of $\mathcal{T}(r)$ or has some some offspring $Q_\ell\in\mathcal{T}(r)$.
    In either case, denote this element by $Q_\ell$.
    Note that all of the $Q_\ell$ are distinct.
    But now since $\delta>0$ is fixed and $\rho(Q)/\rho(Q_0)\geq \delta$, there is a fixed constant $N_\delta$ (independent $m$) so that the cylinder $Q_0$ can be covered by $N_\delta$ balls $B(x,r)$ for $x\in K$.
    By the pigeonhole principle, one of these balls must intersect $m/N_\delta$ cylinders $Q_\ell$.
    Since $m\in\N$ was arbitrary, this contradicts the bounded neighbourhood condition.

    We next bound $\#\mathcal{J}_n$.
    Fix some $Q_0\in\mathcal{T}_n$ be arbitrary and write $r_0 = \rho(Q_0)$, and let $\mathcal{S} = \{Q\in\mathcal{T}: \widehat{Q} = Q_0\}$ denote the set of children of $Q_0$.
    Let $r_0/2 \leq r < r_0$ be sufficiently close to $r_0$ so that $\rho(Q) < r$ for all $Q\in\mathcal{C}$.
    Observe that $\mathcal{S}\subset\mathcal{T}(r)$.
    Since $\diam \pi(Q_0) \leq r_0 \diam X$, there is a constant $N\in\N$ depending only on $\diam X$ so that $\pi(Q)$ can be covered by $N$ balls of radius $r$, say $\{B(x_i, r) \}_{i=1}^N$.
    Since $\pi(Q) \subset \pi(Q_0)$ for all $Q\in\mathcal{S}$, by the bounded neighbourhood condition, each ball $B(x_i, r)$ intersects at most $M$ distinct elements of $\mathcal{S}$.
    Therefore $\#\mathcal{J}_n = \#\mathcal{S} \leq N\cdot M$.
\end{proof}

\subsection{Open set condition and bounds on branching and contraction}
In this section, we discuss the difference between the bounded neighbourhood condition and the open set condition.
Recall the definitions of the open set condition \cref{d:osc} and bounded branching \cref{d:bbranch} from the introduction.

We first give the quick proof relating the bounded branching condition with uniform bounds on $\#\mathcal{J}_n$.
\begin{lemma}\label{l:branch-bound}
    If the non-autonomous IFS $(\Phi_n)_{n=1}^\infty$ satisfies the bounded branching condition, then $\sup_{n\in\N}\#\mathcal{J}_n<\infty$.
    If in addition the IFS is homogeneous, then the converse also holds.
\end{lemma}
\begin{proof}
    In general, if $Q\in\mathcal{T}(r)$ and $Q'$ is a child of $\widehat{Q}$, then $Q'\in\mathcal{T}(r')$ for some $r'\geq r$, so $Q'$ has some offspring in $\mathcal{T}(r)$ which is distinct from $Q$.
    Therefore, $\limsup_{n\to\infty}\#\mathcal{J}_n \le \limsup_{r\to 0}\beta(r)<\infty$.
    Conversely, if the IFS is homogeneous, then for all $r>0$, $\mathcal{T}(r)=\mathcal{T}_{k(r)}$ for some $k(r)\in\N$ so in fact the siblings of $Q$ at scale $r$ are precisely the children of $\widehat{Q}$.
\end{proof}

Finally, we prove \cref{it:reg} which we restate here for convenience.
\begin{restatement}{it:reg}
    Let $(\Phi_n)_{n=1}^\infty$ be a non-autonomous IFS satisfying the open set condition.
    Then the bounded neighbourhood condition holds if any of the following conditions are satisfied:
    \begin{enumerate}[nl,a]
        \item\label{im:osc-bound} $(\Phi_n)_{n=1}^\infty$ satisfies the cone condition and has bounded branching.
        \item\label{im:homog-bound} $(\Phi_n)_{n=1}^\infty$ is homogeneous and has bounded branching.
        \item\label{im:ctr-bound} There is a number $r_{\min}>0$ so that $r_{n,j}\geq r_{\min}$ for all $n\in\N$ and $j\in\mathcal{J}_n$.
    \end{enumerate}
\end{restatement}
\begin{proof}
    Assume that the IFS satisfies the open set condition.
    For notational simplicity, given a cylinder $Q=[j_1,\ldots,j_n]$, we write
    \begin{equation*}
        U(Q) = S_{1,j_1}\circ\cdots\circ S_{n,j_n}(X^\circ).
    \end{equation*}
    Let $x\in K$ and $r>0$ be arbitrary: we need to control the size of $\mathcal{N}(x,r)$.
    We do this under the three assumptions in order.

    First, suppose \cref{im:osc-bound} holds.
    Let $N= \sup_{r>0}\beta(r)$ be the constant from the bounded branching assumption, and moreover observe that it follows from the cone condition that
    \begin{equation*}
        c\coloneqq \inf_{x\in X}\inf_{r\in(0,\diam X)}\frac{\mathcal{L}^d\bigl(B(x,r)\cap X^\circ\bigr)}{\mathcal{L}^d\bigl(B(x,r)\bigr)}>0.
    \end{equation*}
    Write $\mathcal{E}=\{\widehat{Q}:Q\in\mathcal{N}(x,r)\}$.
    Now, let $Q\in\mathcal{N}(x,r)$ be arbitrary and fix some $y_Q\in \pi(Q)\cap B(x,r)$.
    Then since $y_Q\in \pi(\widehat{Q})$ and $\rho(\widehat{Q}) > r$, writing $y_Q = S_{1,j_1}\circ\cdots\circ S_{n-1, j_{n-1}}(z_Q)$ for some $z_Q\in X$ and applying the cone condition to the ball $B(z_Q, r / \rho(\widehat{Q}))$,
    \begin{equation*}
        \mathcal{L}^d\bigl(B(y_Q,r)\cap U(\widehat{Q})\bigr)\geq c\cdot \mathcal{L}^d\bigl(B(x,r)\bigr).
    \end{equation*}
    Moreover, the sets $U(\widehat{Q})\cap B(y_Q,r)\subset B(x, 2r) $ are disjoint for distinct choices of $\widehat{Q}$ by the open set condition, so
    \begin{align*}
        c\cdot\#\mathcal{E}\cdot \mathcal{L}^d\bigl(B(x,r)\bigr) &\leq \sum_{\widehat{Q}\in \mathcal{E}}\mathcal{L}^d\bigl(B(y_Q,r)\cap U(\widehat{Q})\bigr)\\
                                                       &\leq \mathcal{L}^d\bigl(B(x,2r)\bigr)\\
                                                       &= 2^d\mathcal{L}^d\bigl(B(x,r)\bigr).
    \end{align*}
    In other words, $\# \mathcal{E}\leq c^{-1} 2^d$ so by the bounded branching condition, $\#\mathcal{N}(x,r) \leq c^{-1}2^d N$.

    Next, suppose \cref{im:homog-bound} holds.
    Let $N=\sup_{r>0}\beta(r)$ be as before and again set $\mathcal{E}=\{\widehat{Q}:Q\in\mathcal{N}(x,r)\}$.
    Since the IFS is homogeneous, $\rho$ has constant value $r_0$ on $\mathcal{E}$; in particular, if $\widehat{Q}\in \mathcal{E}$, then $\mathcal{L}^d\bigl(U(\widehat{Q})\bigr) = r_0^d \mathcal{L}^d\bigl(X^\circ\bigr)$.
    Moreover, the sets $U(\widehat{Q})$ are disjoint and, since $\pi(\widehat{Q})\cap B(x,r_0)\neq\varnothing$, it follows that $U(\widehat{Q})\subset B(x,2r_0\cdot\diam X)$.
    Therefore
    \begin{align*}
        \#\mathcal{E}\cdot r_0^d\mathcal{L}^d\bigl(X^\circ\bigr) &= \sum_{\widehat{Q}\in \mathcal{E}} \mathcal{L}^d \bigl(U(\widehat{Q})\bigr)\\
                                                                 &\leq \mathcal{L}^d \bigl(B(x,2r_0\cdot \diam X)\bigr)\\
                                                                 & = (2\cdot \diam X)^d r_0^d\mathcal{L}^d(B(0,1)).
    \end{align*}
    Thus $\#\mathcal{N}(x,r)\leq N\cdot (2\cdot\diam X)^d \frac{\mathcal{L}^d(B(0,1))}{\mathcal{L}^d(X^\circ)}$.

    Finally, suppose \cref{im:ctr-bound} holds.
    Now, if $\widehat{Q}\in\mathcal{E}$, since $r r_{\min}\leq \rho(Q)\leq r$, $\pi(Q)\subset B(x,2r\cdot \diam X)$, and the $U(Q)$ are disjoint for distinct $Q$, then
    \begin{align*}
        \#\mathcal{N}(x,r)\cdot (r r_{\min})^d \mathcal{L}^d\bigl(X^\circ\bigr)^d &\leq \sum_{Q\in\mathcal{N}(x,r)}\mathcal{L}^d \bigl(U(Q)\bigr)\\
                                &\leq \mathcal{L}^d\bigl(B(x,2r\cdot\diam X)\bigr)\\
                                & = (2\cdot\diam X)^d r^d\mathcal{L}^d(B(0,1)).
    \end{align*}
    We conclude that $\#\mathcal{N}(x,r)\leq r_{\min}^{-d} (2\diam X)^d \frac{\mathcal{L}^d(B(0,1))}{\mathcal{L}^d(X^\circ)}$.

    Since $x\in K$ and $r>0$ were arbitrary, in any of the cases, it follows that the bounded neighbourhood condition holds.
\end{proof}

\section{Characterizations of Assouad dimension}\label{sec:characterizations}

\subsection{A disc-packing formulation of Assouad dimension}
In this section, we observe that in the definition of the Assouad dimension one may replace the exponent associated to localized coverings of balls of the same size by an exponent coming from localized packings of balls which may have different sizes.
This is essentially the same as the disc-packing formulation for box dimension (see, for example, \cite[§2.6]{zbl:1390.28012}) and our result also proceeds by a similar proof.
This is useful since the natural covers appearing from the symbolic representation of $K$ consist of cylinders which may have very non-uniform diameters when indexed by length.

First, for a metric space $X$, $x\in X$, and $R\in(0,1)$, denote the family of all localized centred packings by
\begin{equation*}
    \operatorname{pack}(X,x,R)=\left\{\{B(x_i,r_i)\}_i:\begin{matrix}0<r_i\leq R,x_i\in X,B(x_i,r_i)\subset B(x,R),\\B(x_i,r_i)\cap B(x_j,r_j)=\varnothing\text{ for all }i\neq j\end{matrix}\right\}.
\end{equation*}
Here the collections may be finite or countably infinite.
We now have the following result.
\begin{proposition}\label{p:Assouad-disc-packing}
    Let $X$ be a bounded metric space.
    Then
    \begin{align*}
        \dimA X=\inf\Bigl\{\alpha:\forall 0<R<1\,\forall x\in X\,&\forall\{B(x_i,r_i)\}_i\in\operatorname{pack}(X,x,R)\\*
                                               &\sum_i r_i^\alpha\lesssim_\alpha R^\alpha\Bigr\}.
    \end{align*}
\end{proposition}
\begin{proof}
    Let $t$ be the right-hand side of the claimed equation.
    That $\dimA X\leq t$ is immediate by specializing to packings with $r_i=r$ for some $0<r\leq R$, using the equivalence (up to a constant factor) of covering and packing counts.

    We now show the lower bound by pigeonholing.
    First, for a centred packing $\mathcal{B}=\{B(x_i,r_i)\}_i\in\operatorname{pack}(X,x,R)$, we write
    \begin{equation*}
        \mathcal{B}_n = \bigl\{B(x_i,r_i)\in\mathcal{B}:2^{-n}\leq r_i/R < 2^{-n+1}\bigr\}.
    \end{equation*}
    Now let $s<t$ and $r_0>0$ be arbitrary.
    Then get some $x\in X$, $0<R<1$, and a centred packing $\mathcal{B}$ as above such that $r_i\leq r_0$ for all $i$ and moreover
    \begin{equation*}
        1\lesssim_s \sum_n \#\mathcal{B}_n 2^{-ns}.
    \end{equation*}
    By the pigeonhole principle with respect to the sequence of weights $(n^{-2})_{n=1}^\infty$, since the weights are summable, there is a $\gamma>0$ depending only on $s$ so that there is some $m\in\N$ with $\gamma\leq \#\mathcal{B}_m m^2 2^{-ms}$.
    Taking $r_0$ sufficiently small, $\#\mathcal{B}_m=0$ for all small $m$, so that $m$ must be arbitrarily large.
    Since for every $\varepsilon>0$ there is a constant $C_\varepsilon$ so that $m^2 \leq C_\varepsilon 2^{m\varepsilon}$ for all $m\in\N$, we conclude that $\dimA K\geq s$.
    Since $s<t$ was arbitrary, we are done.
\end{proof}

\subsection{Subadditivity and submaximality of covering numbers}
In this section, we study a weaker variant of subadditivity, which we call \emph{submaximality}.
Throughout this section, we fix a set $A = G\cap(0,\infty)$ where $G$ is a non-trivial additive subgroup of $\R$.
For concreteness, one might simply fix $A=(0,\infty)$ or $A=\{\kappa n:n\in\N\}$ for some $\kappa>0$.
\begin{definition}\label{d:submax}
    We say that a function $g\colon A\to \R\cup\{-\infty\}$ is \emph{submaximal} if it satisfies the following two assumption:
    \begin{enumerate}[nl,r]
        \item\label{im:subadd-single} For all $y,z\in A$,
            \begin{equation*}
                g(y+z)\leq \max\{g(y),g(z)\}.
            \end{equation*}
        \item\label{im:cont-single} For all $\varepsilon>0$ and $a\in A$, there is an $N=N(\varepsilon, a)$ so that for all $y\geq N$ and $t\leq a$,
            \begin{equation*}
                g(y+t)\leq g(y)+\varepsilon.
            \end{equation*}
    \end{enumerate}
\end{definition}
To motivate this definition, we observe the following version of Fekete's subadditive lemma.
\begin{lemma}\label{l:subadd-single}
    Let $g\colon A\to \R\cup\{-\infty\}$ be submaximal.
    Then
    \begin{equation*}
        \limsup_{y\to\infty}g(y)=\inf_{y\in A}g(y).
    \end{equation*}
\end{lemma}
\begin{proof}
    Let $\varepsilon>0$ and $a\in A$ be arbitrary, and let $N=N(\varepsilon, a)$ be chosen to satisfy the conclusion of \cref{d:submax}~\cref{im:cont-single}.
    Let $y\geq N$ be arbitrary and write $y=\ell a+t$ for $\ell\in\N$ and $0\leq t < a$.
    Then applying \cref{im:cont-single} followed by \cref{im:subadd-single} $\ell-1$ times,
    \begin{equation*}
        g(y) = g(\ell a+t)\leq g(\ell a)+\varepsilon\leq g(a)+\varepsilon.
    \end{equation*}
    Thus
    \begin{equation*}
        \limsup_{y\to\infty}g(y)\leq g(a)+\varepsilon.
    \end{equation*}
    But $a\in A$ and $\varepsilon>0$ were arbitrary, so the desired result follows.
\end{proof}
\begin{remark}
    The assumptions of \cref{l:subadd-single} are satisfied by the function $f(x)/x$, where $f\colon A\to\R\cup\{-\infty\}$ is any subadditive function bounded from above.
    The proof is similar to the proof of \cref{l:subadd-implies-max} below.

    Note that assumption \cref{d:submax}~\cref{im:subadd-single} is not sufficient by itself to guarantee the existence of the limit $\lim_{y\to\infty}g(y)$.
    For example, the function $g\colon\N\to\R$ defined by
    \begin{equation*}
        g(n)=\begin{cases}
            1 &: n\text{ odd},\\
            0 &: n\text{ even},
        \end{cases}
    \end{equation*}
    satisfies $g(y+z)\leq \max\{g(y),g(z)\}$ for all $y,z \in \N$ but has no limit.
\end{remark}

Similarly to before, we say that a function $f\colon A\times A\to\R\cup\{-\infty\}$ is \emph{submaximal} if it satisfies the following two assumptions:
\begin{enumerate}[nl,r]
    \item\label{im:subadd} For all $x,y,z\in A$,
        \begin{equation*}
            f(x,y+z)\leq \max\{f(x,y),f(x+y,z)\}.
        \end{equation*}
    \item\label{im:cont} For all $\varepsilon>0$ and $a\in A$, there is an $N=N(\varepsilon, a)$ so that for all $x,y,t\in A$ with $y\geq N$, $t\leq a$, and $x\leq x'\leq x+t$,
        \begin{equation*}
            f(x, y+t)\leq f(x',y)+\varepsilon.
        \end{equation*}
\end{enumerate}
It is this form that will be essential for us in applications.
We note that the above assumptions are satisfied by the function $f(y,x)/x$, where $f \colon \N \cup \{0\} \times \N \to \R\cup\{-\infty\}$ satisfies the generalized subadditive condition introduced in \cite[§2]{zbl:1078.37014}.
In \cref{l:subadd-implies-max}, we generalize this observation.
But first we acquire the following two-parameter generalization of \cref{l:subadd-single}.
\begin{proposition}\label{l:double-submax}
    Suppose $f\colon A\times A\to\R\cup\{-\infty\}$ is submaximal.
    Then
    \begin{align*}
        \beta\coloneqq{}&\limsup_{y\to\infty}\limsup_{x\to\infty}f(x,y)\\
        ={}&\lim_{y\to\infty}\limsup_{x\to\infty}f(x,y)\\
        ={}&\lim_{y\to\infty}\sup_{x\in A}f(x,y)\\
        ={}&\inf_{y\in A}\sup_{x\in A}f(x,y).
    \end{align*}
    Moreover, if $B\subset A$ is of the form $B=\{\kappa n:n\in\N\}$ for some $\kappa>0$, then
    \begin{equation*}
        \beta=\lim_{\substack{y\to\infty\\y\in B}}\sup_{x\in B}f(x,y).
    \end{equation*}
\end{proposition}
\begin{proof}
    We assume that $\beta>-\infty$: the proof for $\beta=-\infty$ is similar (and substantially easier).
    Write $g(y)=\limsup_{x\to\infty}f(x,y)$.
    We first show that $g$ is submaximal, so that $\beta=\lim_{y\to\infty}g(y)$.
    First,
    \begin{align*}
        g(y_1+y_2)&=\limsup_{x\to\infty}f(x,y_1+y_2)\\
                  &\leq \limsup_{x\to\infty}\max\{f(x,y_1),f(x+y_1,y_2)\}\\
                  &\leq \max\{g(y_1),g(y_2)\}.
    \end{align*}
    Moreover, if $\varepsilon>0$ and $a\in A$ are arbitrary, taking $N$ to satisfy the conclusions of \cref{im:cont} and $y,t\in A$ with $t\leq a$ and $y\geq N$,
    \begin{equation*}
        g(y+t) = \limsup_{x\to\infty}f(x,y+t)\leq \limsup_{x\to\infty}f(x,y)+\varepsilon = g(y)+\varepsilon.
    \end{equation*}
    The same argument with respect to the function $y\mapsto \sup_{x\in A}f(x,y)$ gives that
    \begin{equation*}
        \limsup_{y\to\infty}\sup_{x\in A}f(x,y)=\inf_{y\in A}\sup_{x\in A}f(x,y).
    \end{equation*}

    To complete the proof of the various equalities involving $\beta$, it remains to show that
    \begin{equation}\label{e:beta-eqs}
        \inf_{y\in A}\sup_{x\in A}f(x,y)\leq \beta.
    \end{equation}
    Let $\varepsilon>0$ be arbitrary.
    By definition of $\beta$, there are $y_0$ and $K$ so that for all $x\geq K$, $f(x,y_0)\leq\beta+\varepsilon$.
    Let $N_1=N(\varepsilon, y_0)$ be chosen to satisfy the conclusions of \cref{im:cont} and let $y\in A$ with $y\geq N_1$ be arbitrary.
    Write $y=\ell y_0+t$ for some $\ell\in\N$ and $0\leq t < y_0$.
    Then by \cref{im:subadd}, for all $x\geq K$ and $y\geq N_1$,
    \begin{equation}\label{e:large-bd}
        \begin{aligned}
            f(x,y) &= f(x,\ell y_0+t) \leq  f(x,\ell y_0)+\varepsilon\\
                   &\leq \max_{i=0,\ldots,\ell-1}f(x+i y_0, y_0)+\varepsilon\\
                   &\leq \beta+2\varepsilon.
        \end{aligned}
    \end{equation}
    Now let $N_2=N(\varepsilon, K)$ be chosen to satisfy the conclusions of \cref{im:cont}.
    Let $x\in(0,K)\cap A$ and $y\geq \max\{N_1 + K, N_2\}$.
    Then since $y\geq N_2$ and \cref{e:large-bd} (which applies since $y+x-K\geq N_1$), for all $x\in A$,
    \begin{equation*}
        f(x,y) \leq f(K, y+x-K) + \varepsilon \leq \beta+3\varepsilon.
    \end{equation*}
    Since $\varepsilon>0$ was arbitrary, this proves \cref{e:beta-eqs}.

    Finally, suppose $B\subset A$ is of the form $B=\{\kappa n:n\in\N\}$ for some $\kappa>0$.
    First, note that since $B\subset A$,
    \begin{equation*}
        \beta\geq \lim_{\substack{y\to\infty\\y\in B}}\sup_{x\in B}f(x,y)
    \end{equation*}
    and moreover the limit exists as proven above.
    Conversely, let $(x,y)\in A\times A$ be arbitrary with $y\geq 2\kappa$ and get $(x_0,y_0)\in B\times B$ such that $x\leq x_0<x+\kappa$ and $x+y-\kappa<x_0+y_0\leq x+y$.
    Let $\varepsilon>0$ be arbitrary and let $N=N(\varepsilon,\kappa)$ satisfy the conclusion of \cref{im:cont}.
    Then for $y\geq N$, applying \cref{im:cont} twice,
    \begin{equation*}
        f(x_0,y_0) \geq f(x,y)-2\varepsilon.
    \end{equation*}
    Since $\varepsilon>0$ and $(x,y)\in A\times A$ were arbitrary, it follows that
    \begin{equation*}
        \limsup_{\substack{y_0\to\infty\\y_0\in B}}\sup_{x_0\in B} f(x_0,y_0)\geq \limsup_{y\to\infty}\sup_{x\in A} f(x,y)=\beta
    \end{equation*}
    as required.
\end{proof}
Finally, we show that the hypotheses of \cref{l:double-submax} are satisfied by functions satisfying a two-parameter version of subadditivity.
\begin{lemma}\label{l:subadd-implies-max}
    Suppose $f\colon A\times A\to\{-\infty\}\cup\R$ is any function such that for all $x,y,z\in A$,
    \begin{equation}\label{e:subadd}
        f(x,y+z)\leq \frac{ y f(x,y) + z f(x+y,z)}{y+z}.
    \end{equation}
    Then:
    \begin{enumerate}[nl,r]
        \item\label{im:subadd-impl} For all $x,y,z\in A$,
            \begin{equation*}
                f(x,y+z)\leq \max\{f(x,y), f(x+y,z)\}.
            \end{equation*}
        \item\label{im:cont-impl} Suppose moreover that $f$ is bounded from above by some constant $C>0$.
            Then for all $\varepsilon>0$ and $x,y,t\in A$ with $t\leq \varepsilon C^{-1} y$ and $x\leq x'\leq x+t$,
            \begin{equation*}
                f(x, y+t)\leq f(x',y)+\varepsilon.
            \end{equation*}
            In particular, $f$ is submaximal.
    \end{enumerate}
\end{lemma}
\begin{proof}
    Of course, \cref{im:subadd-impl} is immediate.
    To see \cref{im:cont-impl}, let $C\in\R$ be such that $f(x,y)\leq C$ for all $x,y\in A$.
    Let $\varepsilon>0$.
    Then for all $x,y,t\in A$ with $t\leq \varepsilon C^{-1} y$ and $x\leq x'\leq x+t$, applying \cref{e:subadd} twice,
    \begin{align*}
        f(x,y+t)&\leq \frac{ (x'-x) f(x,x'-x) + (y+t+x-x') f(x',y+t+x-x')}{y+t}\\
                &\leq \frac{ (x'-x) f(x,x'-x) + yf(x',y)+(t+x-x')f(x'+y,t+x-x')}{y+t}\\
                &\leq \frac{y}{y+t} f(x',y)+\frac{t C}{y+t}\\
                &\leq f(x',y) + \varepsilon
    \end{align*}
    as claimed.
\end{proof}
To motivate why the hypotheses of \cref{l:double-submax} and \cref{l:subadd-implies-max} are related to dimensions, we present the observation that the Assouad dimension can be reformulated in a way reminiscent of a notion of dimension studied by Larman \cite{zbl:0152.24502}.
Let $X$ be a bounded doubling metric space and for $\delta\in(0,1)$ and $r\in(0,1)$, write
\begin{equation*}
    \psi(r,\delta) = \sup_{x\in X}N_{r\delta}\bigl(B(x,r)\cap K\bigr)
\end{equation*}
and then set
\begin{equation*}
    \Psi(r,\delta) = \frac{\log \psi(r,\delta)}{\log (1/\delta)}.
\end{equation*}
One can think of $\Psi(r,\delta)$ is the best guess for the Assouad dimension of $X$ at scales $0<r\delta<\delta<1$.
This heuristic is made precise in the following result.
\begin{corollary}\label{c:Assouad-subadd}
    Let $X$ be a bounded doubling metric space.
    Then
    \begin{equation}\label{e:dima-equivs}
        \dimA X=\limsup_{\delta\to 0}\limsup_{r\to 0}\Psi(r,\delta)=\lim_{\delta\to 0}\sup_{r\in(0,1)}\Psi(r,\delta).
    \end{equation}
\end{corollary}
\begin{proof}
    Since $X$ is doubling, there is an $M\geq 0$ so that $\Psi(r,\delta)\in[0,M]$.
    Moreover, given $r,\delta_1,\delta_2\in(0,1)$, by covering balls $B(x,r\delta_1)$ by balls of radius $r\delta_1\delta_2$,
    \begin{equation*}
        \psi(r,\delta_1\delta_2)\leq\psi(r,\delta_1)\psi(r\delta_1,\delta_2)
    \end{equation*}
    and therefore
    \begin{align*}
        \Psi(r,\delta_1\delta_2) &= \frac{\log\psi(r,\delta_1\delta_2)}{\log(1/\delta_1\delta_2)}\\
                                 &\leq\frac{\log\psi(r,\delta_1)+\log\psi(r\delta_1,\delta_2)}{\log(1/\delta_1\delta_2)}\\
                                 &=\frac{\log(1/\delta_1)\Psi(r,\delta_1)+\log(1/\delta_2)\Psi(r\delta_1,\delta_2)}{\log(1/\delta_1)+\log(1/\delta_2)}.
    \end{align*}
    Thus with the change of coordinate $g(x,y)=(e^{-x},e^{-y})$, the second equality in \cref{e:dima-equivs} follows by applying \cref{l:subadd-implies-max} and \cref{l:double-submax} to the function $\Psi\circ g$.

    To see the first equality in \cref{e:dima-equivs}, it is a direct consequence of the definition of the Assouad dimension that
    \begin{equation*}
        \limsup_{\delta\to 0}\limsup_{r\to 0}\Psi(r,\delta)\leq\dimA K
    \end{equation*}
    and that there are sequences $(\delta_n)_{n=1}^\infty$ and $(r_n)_{n=1}^\infty$ with $\lim_{n\to\infty}\delta_n=0$ such that
    \begin{equation*}
        \lim_{\delta\to 0}\sup_{r\in(0,1)}\Psi(r,\delta)\geq\limsup_{n\to\infty}\Psi(r_n,\delta_n)\geq\dimA K,
    \end{equation*}
    as required.
\end{proof}

\section{Assouad dimension of non-autonomous self-similar sets}
In this section, using the results in \cref{sec:reduction-bnc} and \cref{sec:characterizations}, we show that the zero $\theta(n,m)$ of the pressure is submaximal and prove \cref{it:dima}.

\subsection{Regularity of the numbers \texorpdfstring{$\theta(n,m)$}{θ(n,m)}}
Recall the definition of the pressure functions $\phi_{n,m}$ for $n,m\in\N$ in \cref{ss:symb}.
In particular, these functions are convex, differentiable, strictly decreasing, and have unique zero $\theta(n,m)$.
We begin with a upper bound on the derivative of the pressure function at its unique zero.
\begin{lemma}\label{l:pressure-deriv}
    Let $(\Phi_n)_{n=1}^\infty$ be a locally contracting non-autonomous IFS.
    Then there is an $M\in\N$ and a constant $c>0$ so that for all $m\geq M$,
    \begin{equation*}
        \phi_{n,m}'(\theta(n,m)) \leq -c.
    \end{equation*}
\end{lemma}
\begin{proof}
    Since the IFS is locally contracting, there is an $M\in\N$ and a $\delta\in(0,1)$ so that for all $n\in\N$, and $(j_0,\ldots,j_{M-1})\in\mathcal{J}_n\times\cdots\times\mathcal{J}_{n+M-1}$, we have
    \begin{equation*}
        r_{n, j_0}\cdots r_{n+M-1, j_{M-1}} \leq \delta.
    \end{equation*}
    Then for $m\geq M$, write $m = \ell M + j$ for $\ell\in\N$ and we compute
    \begin{align*}
        \phi_{n,m}'(\theta(n,m)) &\leq \frac{1}{m} \sum_{j_0\in\mathcal{J}_n}\cdots \sum_{j_{m-1}\in\mathcal{J}_{n+m-1}}(r_{n, j_0}\cdots r_{n+m-1, j_{m-1}})^{\theta(n,m)}\log \delta^{\ell}\\
                                 &= \log(\delta)\cdot \frac{\ell}{\ell M + j}\\
                                 &\leq - \frac{\log(1/\delta)}{2M}
    \end{align*}
    as claimed.
\end{proof}
Using this, we can now establish submaximality of the function $\theta(n,m)$.
Actually, we will prove a somewhat stronger version of the continuity hypothesis.
\begin{lemma}\label{l:theta-submax}
    Let $(\Phi_n)_{n=1}^\infty$ be a non-autonomous IFS.
    Then for all $n,m,k\in\N$,
    \begin{equation}\label{e:submax}
        \theta(n,m+k) \leq \max\{\theta(n,m),\theta(n+m, k)\}.
    \end{equation}
    Suppose moreover that the IFS satisfies the bounded neighbourhood condition.
    Then there is a constant $C>0$ and an $M\in\N$ so that for all $n,m,k\in\N$ with $m\geq M$ and $n \leq n' \leq n+k$
    \begin{equation*}
        \theta(n, m+k) - \theta(n', m) \leq C\frac{k}{m}.
    \end{equation*}
    In particular, the function $\theta(n,m)$ is submaximal.
\end{lemma}
\begin{proof}
    Let $n,m,k$ be arbitrary.
    Writing $s=\max\{\theta(n,m),\theta(n+m, k)\}$, observe that
    \begin{align*}
        1&=\sum_{j_0\in\mathcal{J}_{n}}\cdots\sum_{j_{m+k-1}\in\mathcal{J}_{n+m+k-1}}\left(\prod_{\ell=0}^{m-1} r_{n+\ell,j_\ell}\right)^{\theta(n,m)}\left(\prod_{\ell=m}^{m+k-1}r_{n+\ell,j_\ell}\right)^{\theta(n+m,\ell)}\\
         &\geq\sum_{j_0\in\mathcal{J}_{n}}\cdots\sum_{j_{m+k-1}\in\mathcal{J}_{n+m+k-1}}\left(\prod_{\ell=0}^{m+k-1} r_{n+\ell,j_\ell}\right)^{s}\\
         &= \exp\bigl((m+k)\cdot\phi_{n,m+k}(s)\bigr).
    \end{align*}
    Since $\phi_{n,m+k}(\theta(n,m+k))=0$ and $\phi_{n,m+k}$ is strictly decreasing, it follows that $\theta(n,m+k)\leq s$, which is \cref{e:submax}.

    We now proceed with the continuity bounds.
    First, by \cref{l:bnc-b}, get $N\in\N$ so that $\#\mathcal{J}_n\leq N$ for all $n\in\N$.
    Now let $n,m\in\N$ be arbitrary.
    To begin, observe that
    \begin{align*}
        \phi_{n,m+1}&(\theta(n,m))\\
                    &=\frac{1}{m+1}\log\sum_{j_0\in\mathcal{J}_{n}}\cdots\sum_{j_{m-1}\in\mathcal{J}_{n+m-1}}\sum_{i\in\mathcal{J}_{n+m}} r_{n,j_0}^{\theta(n,m)}\cdots r_{n+m-1, j_{m-1}}^{\theta(n,m)}\cdot r_{n+m,i}^{\theta(n,m)}\\
                    &= \frac{1}{m+1}\log \sum_{i\in\mathcal{J}_{n+m}} r_{n+m,i}^{\theta(n,m)}\\
                    &\leq \frac{\log N}{m+1}.
    \end{align*}
    Similarly,
    \begin{align*}
        \phi_{n,m+1}(\theta(n+1,m)) &= \frac{1}{m+1}\log \sum_{i\in\mathcal{J}_{n}} r_{n,i}^{\theta(n,m+1)}\leq \frac{\log N}{m+1}.
    \end{align*}
    Next, since the bounded neighbourhood condition implies that the IFS is locally contracting by \cref{l:bnc-b}, applying \cref{l:pressure-deriv}, get a constant $c>0$ and $M\in\N$ so that for all $m\geq M$, $\phi_{n,m+1}'(\theta(n,m+1)) \leq  -c < 0$.
    Since $\phi_{n,m+1}$ is convex, differentiable, and strictly decreasing, since $\phi_{n,m+1} \leq (\log N)/(m+1)$,
    \begin{equation*}
        \min\{\theta(n,m), \theta(n+1,m)\} \geq \theta(n,m+1) - \frac{\log N}{c(m+1)}.
    \end{equation*}
    Applying this 1-step bound $k$ times for general $n,m,k\in\N$ with $m\geq M$ yields the desired result by taking $C = 2(\log N)/c$.
\end{proof}

\subsection{Proof of the Assouad dimension formula}
Finally, we establish our main formula for the Assouad dimension of the non-autonomous self-similar set $K$.
\begin{theorem}\label{t:non-auto-Assouad}
    Let $(\Phi_n)_{n=1}^\infty$ be a non-autonomous IFS satisfying the bounded neighbourhood condition.
    Denote the associated non-autonomous self-similar set by $K$.
    Then
    \begin{equation}\label{e:dimA-theta-form}
        \dimA K = \lim_{m\to\infty}\sup_{n\in\N}\theta(n,m) = \lim_{m\to\infty}\limsup_{n\to\infty}\theta(n,m) = \inf_{m\in\N}\limsup_{n\to\infty}\theta(n,m).
    \end{equation}
\end{theorem}
\begin{proof}
    Recall that the limit in \cref{e:dimA-theta-form} exists by \cref{l:theta-submax} and \cref{l:double-submax}, and moreover
    \begin{equation*}
        s\coloneqq \lim_{m\to\infty}\sup_{n\in\N}\theta(n,m) = \lim_{m\to\infty}\limsup_{n\to\infty}\theta(n,m) = \inf_{m\in\N}\limsup_{n\to\infty}\theta(n,m).
    \end{equation*}
    We verify the lower and upper bounds separately.

    First, recall from \cref{c:symb-Assouad} that $\dimA K=\dimA\Delta$, where $\Delta$ denotes the infinite product space associated with $K$.
    Let $\varepsilon>0$ be fixed and let $M$ be sufficiently large so that for all $m\geq M$, there is an $n\in\N$ so that
    \begin{equation*}
        |\theta(n,m)-s|\leq\varepsilon.
    \end{equation*}
    Now fix a cylinder $[j_1,\ldots,j_n]\subset\Delta$ for some $(j_1,\ldots,j_n)\in\mathcal{J}_1\times\cdots\times\mathcal{J}_n$ and write $R=\diam([j_1,\ldots,j_n])=r_{1,j_1}\cdots r_{n,j_n}$.
    Note that if $m\geq M$, by definition of $\theta(n,m)$
    \begin{equation*}
        \sum_{j_{n+1}\in\mathcal{J}_{n+1}}\cdots\sum_{j_{n+m}\in\mathcal{J}_{n+m}}\prod_{k=1}^{n+m} r_{k,j_k}^{\theta(n,m)}=R^{\theta(n,m)}.
    \end{equation*}
    But the family of cylinders
    \begin{equation*}
        \bigl\{[j_1,\ldots,j_{n+m}]:(j_{n+1},\ldots,j_{n+m})\in\mathcal{J}_{n+1}\times\cdots\times\mathcal{J}_{n+m}\bigr\}
    \end{equation*}
    forms a packing of $B(x,R)$.
    Moreover, since $m\geq M$ is arbitrary, since the IFS is locally contracting by \cref{l:bnc-b}, the width of each cylinder in this family relative to $[j_1,\ldots,j_n]$ converges uniformly to $0$.
    Thus by \cref{p:Assouad-disc-packing}, $\dimA K\geq s-\varepsilon$, so the lower bound follows since $\varepsilon>0$ was arbitrary.

    Conversely, let us upper bound $\dimA K$.
    Let $\varepsilon>0$ be arbitrary, and let $m$ be sufficiently large so that:
    \begin{enumerate}[nl, r]
        \item\label{i:bde} $\theta(n, j) \leq s + \varepsilon$ for all $n\in\N$ and $j \geq m$, and
        \item\label{i:ctr} $r_{n, j_0}\cdots r_{n + m -1, j_{m-1}} \leq 1/2$ for all $n\in\N$ and $j_i \in\mathcal{J}_{n+i}$ for $i=0,\ldots, m-1$.
    \end{enumerate}
    The first choice is possible by definition of $s$, and the second since the IFS is locally contracting by \cref{l:bnc-b}.
    Also by \cref{l:bnc-b}, let $N = \sup_{n\in\N}\#\mathcal{J}_n < \infty$.

    Now let $0<r\leq R<1$ and fix a ball $B(x,R)\subset\Delta$.
    By definition of the metric on $\Delta$, $B(x,R)=[j_1,\ldots,j_p]\eqqcolon Q_0$ where $r_{1,j_1}\cdots r_{p,j_p}\leq R$.
    If $\rho(Q_0) < r$ there is nothing to prove.
    Otherwise, we inductively define $\mathcal{B}_k\subset\mathcal{T}_{p + mk}$ for $k \geq 0$.
    Let $\mathcal{B}_0 = \{Q_0\} \subset \mathcal{T}_p$.
    Now suppose we have defined $\mathcal{B}_k$ for some $k$.
    We then set
    \begin{align*}
        \mathcal{B}_{k+1} = \bigcup_{Q\in\mathcal{B}_k}\{Q'\in\mathcal{T}_{p + m(k+1)}: Q'\subset Q \text{ and }\rho(Q') \geq r\}.
    \end{align*}
    Recall by definition of $\theta$ that for all $k \geq 0$,
    \begin{equation*}
        \sum_{\substack{Q \in \mathcal{T}_{p+m k}\\Q \subset Q_0}}\left(\frac{\rho(Q)}{\rho(Q_0)}\right)^{\theta(p, mk)} = 1.
    \end{equation*}
    Therefore by the choice of $m$ to satisfy \cref{i:bde} and since $r\leq \rho(Q) \leq \rho(Q_0) \leq R$ for all $Q\in\mathcal{B}_k\subset \mathcal{T}_{p+m k}$,
    \begin{equation}\label{e:Bk-bound}
        \#\mathcal{B}_k r^{s+\varepsilon} \leq \sum_{Q\in\mathcal{B}_k} \rho(Q)^{s+\varepsilon} \leq \rho(Q_0)^{s+\varepsilon} \leq R^{s+\varepsilon}.
    \end{equation}
    Let $\ell$ be minimal such that $\mathcal{B}_\ell = \varnothing$; the choice of $m$ to satisfy \cref{i:ctr} implies that $\ell \leq \log_2(R/r) + 1 \leq C(R/r)^\varepsilon$ for some fixed constant $C>0$.

    Next, for $k \geq 0$, let
    \begin{equation*}
        \mathcal{E}_{k+1} = \bigcup_{Q\in\mathcal{B}_k}\{Q'\in\mathcal{T}_{p + m(k+1)}: Q'\subset Q \text{ and }\rho(Q') < r\}.
    \end{equation*}
    This is the set of ($m$-step) children of cylinders in $\mathcal{B}_k$ which are not in $\mathcal{B}_{k+1}$.
    In particular, if for each $Q \in \mathcal{E}_k$ we let $x(Q)\in Q$ (so that $Q\subset B(x(Q), r)$),
    \begin{equation*}
        B(x, R) \subset \bigcup_{k=1}^{\ell}\bigcup_{Q\in\mathcal{E}_k} B(x(Q), r).
    \end{equation*}
    Thus since $\#\mathcal{E}_{k+1} \leq N^m \#\mathcal{B}_{k}$, applying \cref{e:Bk-bound},
    \begin{equation*}
        N_r(B(x,R)) \leq \sum_{k=1}^{\ell}\#\mathcal{E}_k \leq \sum_{k=0}^{\ell - 1} N^m\#\mathcal{B}_{k} \leq CN^m\left(\frac{R}{r}\right)^{s+2\varepsilon}.
    \end{equation*}
    Thus $\dimA K \leq s + 2\varepsilon$.
    Since $\varepsilon>0$ was arbitrary, the upper bound follows.
\end{proof}

\subsection{Some examples exhibiting exceptional behaviour}\label{ss:examples}
In this section, we give two examples of non-autonomous IFSs in $\R$ with respect to the invariant set $X=[0,1]$ with various properties.
In order to lower bound the Assouad dimension in these examples, we find it convenient to use the following lemma.
Let $p_{\mathcal{H}}$ denote the one-sided Hausdorff metric on compact subsets of $\R$:
\begin{equation*}
    p_{\mathcal{H}}(E_1;E_2) = \inf\{\delta>0:\dist(x, E_1) <\delta\text{ for all }x\in E_2\}.
\end{equation*}
The following lemma is standard (see, for instance, \cite[Proposition~6.1.5]{zbl:1201.30002}).
\begin{lemma}\label{l:tan-lower-bd}
    Let $K\subset\R$ be compact, and suppose $E$ is a compact set for which there are numbers $\lambda_n \geq 1$ and $x_n\in\R$ so that
    \begin{equation*}
        \lim_{n\to\infty} p_{\mathcal{H}}(E; \lambda_n(K - x_n)) = 0.
    \end{equation*}
    Then $\dimA K \geq \dimH E$.
\end{lemma}
We now begin with an example demonstrating that the bounded neighbourhood condition is necessary.
First, for $r\in(0,1)$, set
\begin{equation*}
    M(r)=\sup_{x\in K}\#\mathcal{N}(x,r).
\end{equation*}
Then the bounded neighbourhood condition says that $\limsup_{r\to 0}M(r)<\infty$.
\begin{example}\label{ex:unbounded}
    Let $\varepsilon>0$ be arbitrary and let $f(r)$ be any function which diverges to infinity as $r$ converges to $0$.
    We first we give an example of an IFS which:
    \begin{enumerate}[nl,r]
        \item\label{i:osc-cone} satisfies the open set condition and the cone condition,
        \item has $M(r)\leq f(r)$ for all $r$ sufficiently small,
        \item\label{i:size-bd} has $\#\mathcal{J}_n=2$ for all $n\in\N$,
        \item\label{i:dima} has $\dimA K=1$, and
        \item\label{i:small-theta} has $0<\theta(n,m)\leq\varepsilon$ for all $n,m\in\N$.
    \end{enumerate}
    Let us first note some of the implications of the above properties.
    If $\varepsilon<1$, then \cref{t:non-auto-Assouad} together with \cref{i:dima} and \cref{i:small-theta} shows that the bounded neighbourhood condition fails.
    By \cref{it:reg}, this together with \cref{i:osc-cone} implies that the IFS does not have bounded branching.
    Therefore, by \cref{l:branch-bound} and \cref{i:size-bd}, such an IFS cannot be homogeneous.

    Now, we start with the main observation underlying this example, which is to exploit inhomogeneity in a way which makes the numbers $\theta(n,m)$ very small, while still allowing large contraction ratios for flexibility in the construction.
    Given $r_1,r_2\in(0,1)$, set let $s(r_1,r_2)$ denote the unique solution to the equation
    \begin{equation*}
        r_1^{s(r_1,r_2)}+r_2^{s(r_1,r_2)} = 1.
    \end{equation*}
    Then the following is easy to verify: for any $r_1\in(0,1)$,
    \begin{equation}\label{e:small-dim}
        \lim_{r_2\to 0}s(r_1,r_2) = 0.
    \end{equation}

    We now begin the main construction.
    First, fix a number $\ell\in\N$ with $\ell\geq 2$ and a small number $\varepsilon>0$: we define IFSs $\Phi^j_{\ell,\varepsilon}$ for $j=1,\ldots,\ell-1$.
    First, let $r^1_1 = 1 - 1/\ell$, and applying \cref{e:small-dim} let $0<r^1_2< 1-r^1_1$ be sufficiently small so that $s(r^1_1,r^1_2)\leq \varepsilon$, and define $\Phi^1_{\ell,\varepsilon} = \{S^1_1,S^1_2\}$ by
    \begin{equation*}
        S^1_1(x) = r^1_1 x\qquad\text{and}\qquad S^1_2(x) = r^1_2 x + r^1_1.
    \end{equation*}
    Now suppose we have defined $\Phi^j_{\ell,\varepsilon}$ for some $j\in\N$ with $j<\ell-1$.
    Let $r^{j+1}_1$ be chosen so that $r^1_1\cdots r^j_1\cdot r^{j+1}_1 = 1 - \frac{j+1}{\ell}$, and let $r^{j+1}_2 < 1-r^{j+1}_1$ be chosen so that $s(r^{j+1}_1, r^{j+1}_2) \leq\varepsilon$, and define $\Phi^{j+1}_{\ell,\varepsilon} = \{S^{j+1}_1,S^{j+1}_2\}$ by
    \begin{equation*}
        S^{j+1}_1(x) = r^{j+1}_1 x\qquad\text{and}\qquad S^{j+1}_2(x) = r^{j+1}_2 x + r^{j+1}_1.
    \end{equation*}
    Note for each $j=1,\ldots,\ell-1$ that
    \begin{equation}\label{e:2-choice}
        S^1_1\circ\cdots\circ S^{j-1}_1\circ S^j_2(0) = r^1_1\cdots r^{j-1}_1 r^j_1 = 1 - \frac{j}{\ell}.
    \end{equation}
    Observe that for every $j=1,\ldots,\ell-1$ the IFS $\Phi^j_{\ell,\varepsilon} = \{S^j_1, S^j_2\}$ we have defined satisfies the open set condition with respect to the interval $(0,1)$ since $r^j_1+r^j_2 <1$ and the similarity dimension satisfies $s(r^{j}_1, r^{j}_2) \leq\varepsilon$.

    Next, for each $\ell\in\N$ with $\ell\geq 2$, define the tuple of IFSs
    \begin{equation*}
        \Phi^\ell = \bigl(\Phi^1_{\ell, \varepsilon},\ldots,\Phi^{\ell-1}_{\ell, \varepsilon}\bigr).
    \end{equation*}
    Our non-autonomous IFS will be a concatenation of the sequence of tuples of IFSs $(\Phi^{\ell_k})_{k=1}^\infty$, for numbers $\ell_k$ with $\limsup_{k\to\infty}\ell_k = \infty$.
    We will choose these numbers in blocks $(\ell_1,\ldots,\ell_{m_k})$ for $m_k\in\N$.
    To begin, the IFS $\Phi^1_{2,\varepsilon}$ satisfies the bounded neighbourhood condition with some constant $N_1$.
    Let $r_1$ be such that $f(r) \geq 2^2 \cdot N_1$ for all $0<r\leq r_1$, and choose $m_1$ sufficiently large so that with
    \begin{equation*}
        (\ell_1,\ldots,\ell_{m_1}) = (2,\ldots,2)
    \end{equation*}
    every cylinder corresponding to a composition from $(\Phi^{\ell_1},\ldots,\Phi^{\ell_{m-1}})$ has contraction ratio at most $r_1$.

    Now suppose we have chosen $(\ell_1,\ldots,\ell_{m_k})$ for some $k\in\N$.
    By \cref{it:reg}, the non-autonomous IFS formed by repeating the maps in $(\Phi^{\ell_1},\ldots,\Phi^{\ell_{m_k}})$ satisfies the bounded neighbourhood condition with some constant $N_k$.
    Since $f(r)$ diverges to infinity, get $0<r_{k+1}\leq r_k$ so that $f(r)\geq N_k 2^{k+1}$ for all $0<r\leq r_{k+1}$.
    Then let $n_k$ be sufficiently large so that every cylinder corresponding to a composition in $(\ell_1,\ldots,\ell_{m_k})^{n_k}$ has contraction ratio at most $r_{k+1}$.
    Then let $m_{k+1} = n_k m_k + 1$ and define
    \begin{equation*}
        (\ell_1,\ldots,\ell_{m_{k+1}}) = (\ell_1,\ldots,\ell_{m_k},\ldots,\ell_1,\ldots,\ell_{m_k}, k+2).
    \end{equation*}

    Let $\Phi$ denote the infinite concatenation of $(\Phi^{\ell_k})_{k=1}^\infty$ have corresponding function $M(r)$.
    For each $k\in\N$, since the IFS corresponding to the tuple $\Phi^{k+2}$ has $2^{k+1}$ maps, it follows that $N_{k+1} \leq N_k 2^{k+1}$.
    Moreover, by construction, for all $0<r_{k+1}<r \leq r_k$, $M(r) \leq N_k$; and on the other hand, the $r_k$ were chosen precisely so that $f(r) \geq N_{k-1} 2^k \geq N_k$ for all $0<r\leq r_k$.
    Therefore, $M(r) \leq f(r)$ for all $0<r\leq r_1$, as required.

    Moreover, by construction, this non-autonomous IFS satisfies the open set condition and $\#\mathcal{J}_n=2$ for all $n$.
    Let $K$ denote the corresponding limit set.

    Observe that $0\in K$ since $0$ is a fixed point of the first map in each component IFS.
    Next, for $k\in\N$, let $g$ denote the map corresponding to any cylinder from $(\Phi^{\ell_1},\ldots,\Phi^{\ell_{m_k-1}})$.
    Since $0\in K$, by \cref{e:2-choice},
    \begin{equation*}
        \left\{\frac{1}{k+1},\ldots,\frac{k}{k+1}\right\}\subset g^{-1}(K)\cap [0,1].
    \end{equation*}
    Since $k\in\N$ was arbitrary, by \cref{l:tan-lower-bd}, $\dimA K=1$.

    Finally, to verify the condition concerning the numbers $\theta(n,m)$, by construction $\theta(n,1)\leq\varepsilon$ for all $n\in\N$, so by \cref{l:theta-submax}, $\theta(n,m)\leq\varepsilon$ for all $n,m$.
\end{example}
Next, we give an example which is homogeneous and satisfies the open set condition and the cone condition for which the Assouad dimension and the symbolic formula for the Assouad dimension can be specified arbitrarily.
\begin{example}\label{ex:arbitrary-values}
    Let $(k_n)_{n=1}^\infty$ be an arbitrary sequence with $\limsup_{n\to\infty}k_n=\infty$.
    Let $0<s\leq t \leq 1$.
    We give an example of an IFS which:
    \begin{enumerate}[nl,r]
        \item is homogeneous,
        \item satisfies the open set condition,
        \item has $\#\mathcal{J}_n\leq k_n$ for all $n\in\N$,
        \item\label{i:3} has $\dimA K=t$, and
        \item\label{i:4} has $\theta(n,m)=s$ for all $n,m\in\N$.
    \end{enumerate}
    Observe that if $s<t$, then \cref{t:non-auto-Assouad} together with \cref{i:3} and \cref{i:4} implies that the bounded neighbourhood condition fails.
    Alternatively, if $\limsup_{n\to\infty}\#\mathcal{J}_n=\infty$, then \cref{l:bnc-b} implies that the bounded neighbourhood condition fails.
    The idea in this construction is to successively approximate a Cantor set of dimension $t$ at large levels $n$, using maps with ratios corresponding to a Cantor set of dimension $s$.
    Let $r_1$ be chosen so that $\log 2 / \log(1/r_1) = s$ and let $r_2$ be chosen so that $\log 2/ \log(1/r_2)=t$.
    Also, let $\Psi = \{S_1,S_2\}$ denote the Cantor IFS where $S_1(x) = r_2 x$ and $S_2(x) = r_2 x + (1-r_2)$, and let $E_n$ denote the set of left endpoints at level $n$:
    \begin{equation*}
        E_n = \bigl\{S_{i_1}\circ\cdots\circ S_{i_n}(0):(i_1,\ldots,i_n)\in\{1,2\}^n\bigr\}
    \end{equation*}
    We then denote corresponding maps
    \begin{equation*}
        \Phi_n = \{ x\mapsto r_1^n x + y: y\in E_n\}.
    \end{equation*}
    Finally, let $(m_n)_{n=1}^\infty$ be a sequence of natural numbers with $\limsup_{n\to\infty}m_n = \infty$, and $2^{m_n} \leq k_n$ for all $n\in\N$.
    Such a choice is possible since, by assumption, $\limsup_{n\to\infty}k_n = \infty$.
    We then define the non-autonomous IFS
    \begin{equation*}
        \Phi = (\Phi_{m_n})_{n=1}^\infty.
    \end{equation*}
    Since $\#\mathcal{J}_n = 2^{m_n} \leq k_n$, the growth rate condition is satisfied, and the IFS is also clearly homogeneous since all of the contraction ratios in $\Phi_{m_n}$ are exactly $r_1^{m_n}$.
    Moreover, observe that $\theta(n,1)=s$ by choice of $r_1$ for all $n\in\N$, so in fact $\theta(n,m) = s$ for all $n,m\in\N$.

    Let $K$ denote the corresponding limit set: it remains to show that $\dimA K=t$.
    First note that $0\in K$, since $0\in E_n$ for all $n\in\N$.
    Now let $M\in\N$ be arbitrary and get $n$ such that $m_n \geq M$.
    Let $f= f_1\circ\cdots\circ f_{n-1}$ for $f_j \in \Psi_j$.
    Since $0\in K$, $E_M\subset E_{m_n} \subset f^{-1}(K)\cap [0,1]$.
    But $\lim_{M\to\infty} E_M = E$ in the Hausdorff metric, where $E$ is the attractor of the IFS $\Psi$.
    Therefore by \cref{l:tan-lower-bd}, $\dimA K \geq \dimH E = t$.

    Now, we obtain the upper bound.
    Write $\rho_n = r_1^{m_1+\cdots+m_n}$.
    We use a different covering strategy between different pairs of scales.
    \begin{enumerate}
        \item We first observe between pairs of scales $\rho_n$ and $\rho_{n+k}$ that $K$ has dimension approximately $s$.
            More precisely, let $f\in\Phi_1\circ\cdots\circ\Phi_n$ and consider the interval $I = f([0,1])$.
            Then $I$ can be covered by $2^{m_{n+1}+\cdots+m_{n+k}}$ intervals of side-length $\rho_{n+k}$, so that
            \begin{equation*}
                \frac{\log N_{\rho_{n+k}}(K\cap I)}{\log(\rho_n/\rho_{n+k})} \leq \frac{(m_{n+1}+\cdots m_{n+k})\log 2}{(m_{n+1}+\cdots m_{n+k})\log(1/r_1)} = s.
            \end{equation*}
            But any ball $B(x,\rho_n)$ can be covered by at most $2$ intervals $I=f([0,1])$, so
            \begin{equation}\label{e:outer-scales}
                N_{\rho_{n+k}}(K\cap B(x, \rho_n)) \leq 2\left(\frac{\rho_n}{\rho_{n+k}}\right)^s.
            \end{equation}
        \item Next, let $\rho_n \geq R\geq r \geq \rho_{n+1}$ and consider some interval $I = f([0,1])$ where $f\in\Phi_1\circ\cdots\circ\Phi_n$.
            Since $E$ is a self-similar set satisfying the strong separation condition with dimension $t$, there is a constant $C_0>0$ so that
            \begin{equation*}
                N_{r/\rho_n}\bigl(E \cap B(x,R/\rho_n)\bigr) \leq C_0\left(\frac{R}{r}\right)^t.
            \end{equation*}
            Observe moreover that $K\cap f([0,1])$ is contained in a $\rho_{n+1}$-neighbourhood of $f(E_{m_n})\subset f(E)$.
            But $B(x,R)$ is contained in at most $2$ intervals of the form $I = f([0,1])$, so there is a constant $C>0$ so that
            \begin{equation}\label{e:inner-scales}
                N_r(B(x,R)\cap K) \leq C\left(\frac{R}{r}\right)^t.
            \end{equation}
    \end{enumerate}
    Let $0<r\leq R<1$ be arbitrary and let $x\in K$.
    Let $n_1\in\N$ be maximal so that $\rho_{n_1}\geq R$ and let $n_2\in\N$ be maximal so that $\rho_{n_2}\geq r$.
    If there is an $n$ such that $\rho_n \geq R\geq r \geq\rho_{n+1}$, we are simply done by \cref{e:outer-scales}.
    Otherwise, there is a minimal $n$ and a maximal $k$ so that $R \geq \rho_n \geq \rho_{n+k} \geq r$.
    Applying \cref{e:inner-scales} to the pairs $R\geq \rho_n$ and $\rho_{n+k}\geq r$, and \cref{e:outer-scales} between $\rho_n\geq \rho_{n+k}$ yields the bound
    \begin{align*}
        N_r\bigl(B(x,R)\cap K\bigr) &\leq C \left(\frac{R}{\rho_n}\right)^t \cdot 2\left(\frac{\rho_n}{\rho_{n+k}}\right)^s \cdot C \left(\frac{\rho_{n+k}}{r}\right)^s\\
                                    & \leq 2 C^2 \left(\frac{R}{r}\right)^t.
    \end{align*}
    Therefore $\dimA K\leq t$, as claimed.
\end{example}
\begin{acknowledgements}
    The authors thank Balázs Bárány and Lars Olsen for many valuable comments on a draft version of this document, in particular for pointing out an error in the proof of an earlier version of \cref{t:non-auto-Assouad}.

    AR is supported by Tuomas Orponen's grant from the Research Council of Finland via the project Approximate incidence geometry, grant no.\ 355453.
    Part of this work was also done while AR was a PhD student at the University of St Andrews, supported by EPSRC Grant EP/V520123/1 and the Natural Sciences and Engineering Research Council of Canada.

    An earlier version of \cref{l:bnc-b} claimed that the bounded neighbourhood condition implies bounded branching.
    This is false.
    We thank Minghui Xu for pointing out the error in the proof.
\end{acknowledgements}
\end{document}